%% file: Eigenvector_monotonicity_v3.tex
\documentclass[12pt]{article}
\usepackage[latin1]{inputenc}
\usepackage[british]{babel}
\usepackage{cmap}
\usepackage{lmodern}

\usepackage{amssymb, amsmath, amsthm}
\usepackage[a4paper,top=25mm,bottom=25mm,left=25mm,right=25mm]{geometry}
\usepackage{etex}
\usepackage{ragged2e}

\usepackage{authblk} 
\usepackage{pifont}
\usepackage{graphicx}
\usepackage[usenames,dvipsnames,svgnames,table]{xcolor}
\usepackage[figuresright]{rotating}
\usepackage{xtab} 
\usepackage{longtable} 
\usepackage{multirow}
\usepackage{footnote}
\usepackage[stable]{footmisc}
\usepackage{chngpage} 
\usepackage{pdflscape} 
\usepackage[nottoc,notlot,notlof]{tocbibind} 

\usepackage{pgfplots}
\pgfplotsset{every tick label/.append style={font=\footnotesize}}
\pgfplotsset{compat=1.14}
\usepackage{setspace}

\usepackage{array}
\newcolumntype{K}[1]{>{\centering\arraybackslash$}p{#1}<{$}}

\makesavenoteenv{tabular}
\usepackage{tabularx}
\usepackage{booktabs}
\usepackage{threeparttable}
\usepackage[referable]{threeparttablex} 
\newcolumntype{R}{>{\raggedleft\arraybackslash}X}
\newcolumntype{L}{>{\raggedright\arraybackslash}X}
\newcolumntype{C}{>{\centering\arraybackslash}X}
\newcolumntype{A}{>{\columncolor{gray!25}}C}
\newcolumntype{a}{>{\columncolor{gray!25}}c}

\newlength{\tablen}

\usepackage{dcolumn} 
\newcolumntype{.}{D{.}{.}{-1}}

\usepackage{tikz}
\usetikzlibrary{arrows, calc, matrix, patterns, positioning, trees}
\usepackage[semicolon]{natbib}
\usepackage[hyphens]{url}
\usepackage{hyperref} 
\hypersetup{
  colorlinks   = true,    		
  urlcolor     = blue,    		
  linkcolor    = blue,    		
  citecolor    = ForestGreen	
}
\usepackage{microtype}
\usepackage[justification=centering]{caption} 

\usepackage[labelformat=simple]{subcaption}

\DeclareCaptionLabelFormat{parenthesis}{(#2)}
\makeatletter
\renewcommand\p@subfigure{\arabic{figure}.}
\makeatother

\DeclareCaptionLabelFormat{parenthesis}{(#2)}
\makeatletter
\renewcommand\p@subtable{\arabic{table}.}
\makeatother

\usepackage{enumitem}

\setlist[itemize]{leftmargin=2.5\parindent}
\setlist[enumerate]{leftmargin=2.5\parindent}

\newenvironment{customlegend}[1][]{%
	\begingroup
	\csname pgfplots@init@cleared@structures\endcsname
	\pgfplotsset{#1}%
    }{%
	\csname pgfplots@createlegend\endcsname
	\endgroup
    }%
\def\addlegendimage{\csname pgfplots@addlegendimage\endcsname}

\theoremstyle{plain}

\newtheorem{proposition}{Proposition}[section]

\theoremstyle{definition}
\newtheorem{axiom}{Axiom}

\newtheorem{definition}{Definition}[section]
\newtheorem{example}{Example}[section]

\theoremstyle{remark}

\newtheorem{remark}{Remark}

\makeatletter
\let\@fnsymbol\@alph
\makeatother

\def\keywords{\vspace{.5em} 
{\noindent \textit{Keywords}: }}

\def\AMS{\vspace{.5em} 
{\noindent \textbf{\emph{MSC} class}: }}

\def\JEL{\vspace{.5em} 
{\noindent \textbf{\emph{JEL} classification number}: }}

\title{On the monotonicity of the eigenvector method}
\author{\href{https://sites.google.com/view/laszlocsato}{L\'aszl\'o Csat\'o}\thanks{~Corresponding author. E-mail: \emph{laszlo.csato@sztaki.hu} \newline
Institute for Computer Science and Control (SZTAKI), Laboratory on Engineering and Management Intelligence, Research Group of Operations Research and Decision Systems \newline
Corvinus University of Budapest (BCE), Department of Operations Research and Actuarial Sciences} $\qquad \qquad$
\href{https://sites.google.com/view/doragretapetroczy}{D\'ora Gr\'eta Petr\'oczy}\thanks{~E-mail: \emph{doragreta.petroczy@uni-corvinus.hu} \newline Corvinus University of Budapest (BCE), Department of Finance}} 

\date{\today}

\def\Dedication{
\begin{small}
{\noindent
``\emph{But nobody should deny the following principle}''\footnote{~In the original: ``\emph{Niemand aber d\"urfte folgendes Prinzip bestreiten}'' \citep[p.~201]{Landau1914}.}
}
\end{small}

\flushright
\begin{small}
(Edmund Landau: \emph{\"Uber Preisverteilung bei Spielturnieren})
\end{small}

\vspace{0.5cm} 
\justify }

\begin{document}

\maketitle
\thispagestyle{empty}
\Dedication

\begin{abstract}
\noindent
Pairwise comparisons are used in a wide variety of decision situations where the importance of alternatives should be measured on a numerical scale. One popular method to derive the priorities is based on the right eigenvector of a multiplicative pairwise comparison matrix.
We consider two monotonicity axioms in this setting. First, increasing an arbitrary entry of a pairwise comparison matrix is not allowed to result in a counter-intuitive rank reversal, that is, the favoured alternative in the corresponding row cannot be ranked lower than any other alternative if this was not the case before the change (rank monotonicity). Second, the same modification should not decrease the normalised weight of the favoured alternative (weight monotonicity).
Both properties are satisfied by the geometric mean method but violated by the eigenvector method. The axioms do not uniquely determine the geometric mean. 
The relationship between the two monotonicity properties and the Saaty inconsistency index are investigated for the eigenvector method via simulations. Even though their violation turns out not to be a usual problem even for heavily inconsistent matrices, all decision-makers should be informed about the possible occurrence of such unexpected consequences of increasing a matrix entry.


\AMS{90B50, 91B08}

\JEL{C44, D71}
\end{abstract}

\section{Introduction} \label{Sec1}

Several decision-making methods involve the comparison of the criteria and the alternatives in pairs, making judgements, and compiling the results into multiplicative positive reciprocal pairwise comparison matrices. For instance, this is a crucial component of the Analytic Hierarchy Process (AHP), introduced by Thomas L. Saaty \citep{Saaty1977, Saaty1980}. He has suggested deriving the priorities from such a matrix by its principal right eigenvector, which is called the \emph{eigenvector method}.
Since AHP has numerous applications \citep{Ho2008, SaatyVargas2012, VaidyaKumar2006}, a better understanding of this procedure seems to be a prominent research question.

The starting point of our paper is a remark of Saaty \citep[p.~86]{Saaty2003}: ``\emph{Now we ask the question, what is priority or more generally what meaning should we attach to a priority vector of a set of alternatives? We can think of two meanings. The first is a numerical ranking of the alternatives that indicates an order of preference among them. The other is that the ordering should also reflect intensity or cardinal preference as indicated by the ratios of the numerical values and is thus unique to within a positive multiplicative constant (a similarity transformation).}''
A potential implication is that both the rank and the (normalised) weight of any alternative should be a \emph{monotonic} function of its comparisons.
In other words, increasing an arbitrary entry in the $i$th row of a pairwise comparison matrix should not result in a rank reversal such that alternative $i$ was ranked at least as high as alternative $k$ before the change but it is ranked lower after the change. Similarly, such a change is not allowed to decrease the overall importance of alternative $i$.
If the straightforward conditions of \emph{rank monotonicity} or \emph{weight monotonicity} do not hold, then the ordering or the weights of the alternatives may behave against the intentions of the decision-maker who wants to express a stronger preference for an alternative by increasing its pairwise comparisons.

The eigenvector method will be proved to suffer from both problems for certain pairwise comparison matrices, while the geometric mean method always satisfies these requirements. It is also shown that the axioms do not uniquely determine the geometric mean, hence they can serve as reasonable criteria to classify the weighting methods suggested for pairwise comparison matrices.

It will also be investigated how often the eigenvector method violates these properties for a randomly generated matrix as a function of its consistency ratio, the consistency measure suggested by \citet{Saaty1977}. In particular, the axioms are tested by substituting a randomly chosen matrix entry with the next element on the well-known Saaty scale.
Their violation turns out not to be a usual problem even for heavily inconsistent matrices.
Nonetheless, all decision-makers should be informed about the possible occurrence of counter-intuitive changes after a matrix entry is increased because probably nobody expects that a more favourable opinion on an alternative can be detrimental to its rank or weight.

The paper has the following structure.
Section~\ref{Sec2} outlines the topic of pairwise comparison matrices and introduces the two axioms of monotonicity. The eigenvector method and some other procedures are analysed in the view of these properties in Section~\ref{Sec3}.
Finally, Section~\ref{Sec4} concludes.

\section{The problem} \label{Sec2}

In this section, the main notions around pairwise comparison matrices are briefly recalled, and two natural properties are presented. We also offer a short overview of related papers.

\subsection{Preliminaries: multiplicative pairwise comparison matrices} \label{Sec21}

Let $N = \{ 1,2, \dots ,n \}$ be a set of alternatives to be evaluated. Assume that their pairwise comparisons are known: $a_{ij}$ is a numerical answer to the question ``\emph{How many times alternative $i$ is better than alternative $j$?}'', that is, $a_{ij}$ quantifies the relative importance of alternative $i$ with respect to alternative $j$.

Let $\mathbb{R}^{n}_+$ and $\mathbb{R}^{n \times n}_+$ denote the set of positive (with all elements greater than zero) vectors of size $n$ and matrices of size $n \times n$, respectively.

The comparisons are collected into a matrix whose entries below the diagonal are reciprocal to the corresponding entries above the diagonal.

\begin{definition} \label{Def21}
\emph{Multiplicative pairwise comparison matrix}:
Matrix $\mathbf{A} = \left[ a_{ij} \right] \in \mathbb{R}^{n \times n}_+$ is a \emph{multiplicative pairwise comparison matrix} if $a_{ji} = 1/a_{ij}$ for all $1 \leq i,j \leq n$.
\end{definition}

In the following, the word ``multiplicative'' will be omitted for the sake of simplicity.

The set of all pairwise comparison matrices with $n$ alternatives is denoted by $\mathcal{A}^{n \times n}$.

Pairwise comparisons are carried out in order to obtain a priority vector $\mathbf{w}$ such that the proportion of the weights $w_i$ and $w_j$ of the alternatives $i$ and $j$, respectively, approximates the value of their pairwise comparison, that is, $w_i / w_j \approx a_{ij}$. Thus the weights can be normalised arbitrarily.

\begin{definition} \label{Def22}
\emph{Weight vector}:
Vector $\mathbf{w}  = \left[ w_{i} \right] \in \mathbb{R}^n_+$ is a \emph{weight vector} if $\sum_{i=1}^n w_{i} = 1$.
\end{definition}

The set of weight vectors of size $n$ is denoted by $\mathcal{R}^{n}$.

\begin{definition} \label{Def23}
\emph{Weighting method}:
Function $f: \mathcal{A}^{n \times n} \to \mathcal{R}^{n}$ is a \emph{weighting method}.
\end{definition}

The weight of alternative $i$ from the pairwise comparison matrix $\mathbf{A}$ according to the weighting method $f$ is denoted by $f_i(\mathbf{A})$.

Weighting methods are often used to rank the alternatives.
Ranking $\succeq$ is a weak order on the set of alternatives $N$.
Any ranking $\succeq$ has two parts, the asymmetric relation $\succ$, and the symmetric relation $\sim$, defined as follows: $i \succ j$ if and only if $i \succeq j$ but $i \preceq j$ does not hold, and $i \sim j$ if and only if $i \succeq j$ and $i \preceq j$, respectively.

There exist many methods to estimate a suitable weight vector from a pairwise comparison matrix. Probably the most popular procedures are the \emph{(row) geometric mean (logarithmic least squares) method} \citep{WilliamsCrawford1980, CrawfordWilliams1985, DeGraan1980, deJong1984, Rabinowitz1976}, and the eigenvector method \citep{Saaty1977}. Although the latter suffers from a number of theoretical shortcomings discussed later, and there are sound axiomatic arguments in favour of the geometric mean \citep{Fichtner1984, LundySirajGreco2017, Csato2018c, BozokiTsyganok2019, Csato2019a}, the AHP methodology mainly uses the eigenvector method since the pioneering work of Saaty. Therefore, this procedure will be in our focus.

\begin{definition} \label{Def24}
\emph{Eigenvector method} \citep{Saaty1977}:
The \emph{eigenvector method} associates the weight vector $\mathbf{w}^{EM} (\mathbf{A}) \in \mathcal{R}^n$ for a given pairwise comparison matrix $\mathbf{A} \in \mathcal{A}^{n \times n}$ such that
\begin{equation} \label{eq_EM}
\mathbf{A} \mathbf{w}^{EM}(\mathbf{A}) = \lambda_{\max}(\mathbf{A}) \mathbf{w}^{EM}(\mathbf{A}),
\end{equation}
where $\lambda_{\max}(\mathbf{A})$ denotes the maximal eigenvalue, also known as the principal or Perron eigenvalue, of the (positive) matrix $\mathbf{A}$.
\end{definition}

\begin{definition} \label{Def25}
\emph{Geometric mean method} \citep{WilliamsCrawford1980, CrawfordWilliams1985, DeGraan1980, deJong1984, Rabinowitz1976}:
The \emph{geometric mean method} associates the weight vector $\mathbf{w}^{GM} (\mathbf{A}) \in \mathcal{R}^n$ for a given pairwise comparison matrix $\mathbf{A} \in \mathcal{A}^{n \times n}$, where
\begin{equation} \label{eq_GM}
w_i^{GM}(\mathbf{A}) = \frac{\prod_{j=1}^n a_{ij}^{1/n}}{\sum_{k=1}^n \prod_{j=1}^n a_{kj}^{1/n}}.
\end{equation}
\end{definition}

Another priority deriving method is considered to illustrate monotonicity.

\begin{definition} \label{Def26}
\emph{Column sum method} \citep{Zahedi1986, ChooWedley2004}:
The \emph{column sum method} associates the weight vector $\mathbf{w}^{CM} (\mathbf{A}) \in \mathcal{R}^n$ for a given pairwise comparison matrix $\mathbf{A} \in \mathcal{A}^{n \times n}$, where
\begin{equation} \label{eq_CM}
w_i^{CM}(\mathbf{A}) = \frac{\sum_{j=1}^n a_{ij}}{\sum_{k=1}^n \sum_{j=1}^n a_{kj}}.
\end{equation}
\end{definition}

All weighting methods induce a ranking, for instance, $i \succeq_{\mathbf{A}}^{EM} j$ if and only if $w_i^{EM}(\mathbf{A}) \geq w_j^{EM}(\mathbf{A})$.

There is a special case when all reasonable weighting methods, including the above three, give the same result.

\begin{definition} \label{Def28}
\emph{Consistency}:
Let $\mathbf{A} = \left[ a_{ij} \right] \in \mathbb{R}^{n \times n}_+$ be a pairwise comparison matrix. It is called \emph{consistent} if the condition $a_{ik} = a_{ij} a_{jk}$ holds for all $1 \leq i,j,k \leq n$.
\end{definition}

However, consistency is seldom observed in practice, pairwise comparison matrices are usually \emph{inconsistent}. A variety of indices has been proposed to measure the level of inconsistency, see \citet{Brunelli2018} for a survey of them.
We will consider the oldest and by far the most popular Saaty inconsistency index \citep{Saaty1977}, which is closely related to the eigenvector method.

\begin{definition} \label{Def29}
\emph{Consistency index} ($CI$):
Let $\mathbf{A} = \left[ a_{ij} \right] \in \mathbb{R}^{n \times n}_+$ be a pairwise comparison matrix. Its \emph{consistency index} is
\[
CI(\mathbf{A}) = \frac{\lambda_{\max}(\mathbf{A})-n}{n-1},
\]
where $\lambda_{\max}(\mathbf{A})$ is the principal eigenvalue of matrix $\mathbf{A}$ as before.
\end{definition}

\citet{Saaty1977} introduced the so-called \emph{random index} $RI_n$, that is, the average $CI$ of a large number of $n \times n$ pairwise comparison matrices with entries randomly generated from the scale $\{ 1/9, 1/8, \dots ,8,9 \}$. The proportion of $CI$ and $RI_n$ is called the \emph{consistency ratio} $CR$, or the Saaty inconsistency index.

\citet{Saaty1977} considered a pairwise comparison matrix to be acceptable if the value of $CR$ does not exceed the threshold $0.1$.

\subsection{Monotonicity on single comparisons} \label{Sec22}

The entry $a_{ij}$ measures the dominance of alternative $i$ over alternative $j$, thus it is not expected that increasing $a_{ij}$ leads to an ordering where alternative $i$ is ranked lower than any alternative $k$ if it was ranked at least as high before the change. Such counter-intuitive rank reversal would be probably against the intention of the decision-maker.

The following axiom formalises this requirement.

\begin{axiom} \label{Axiom1}
\emph{Rank monotonicity}:
Let $\mathbf{A} \in \mathcal{A}^{n \times n}$ be any pairwise comparison matrix and $1 \leq i,j \leq n$ be any two different alternatives. Let $\mathbf{A'} \in \mathcal{A}^{n \times n}$ be identical to $\mathbf{A}$ but $a_{ij}' > a_{ij}$ (and $a_{ji}' < a_{ji}$ due to the reciprocity property). 
The weighting method $f: \mathcal{A}^{n \times n} \to \mathcal{R}^n$ is called \emph{rank monotonic} if $i \succeq_{\mathbf{A}}^f k \Rightarrow i \succeq_{\mathbf{A'}}^f k$ for all $1 \leq k \leq n$.
\end{axiom}

Rank monotonicity is widely considered in social choice theory, hence it would be difficult to give credit to any author for first suggesting this property. \citet{Rubinstein1980} uses it to characterise the ranking method commonly known as the points system on the set of tournaments, where either player $i$ beats player $j$, or player $j$ beats player $i$. \citet{Henriet1985} extends this result by dropping the asymmetric property of the beating relation, which is further generalised by \citet{Bouyssou1992} to the case of valued relations. The axiom is called \emph{positive responsiveness} in \citet{vandenBrinkGilles2009} and \emph{positive responsiveness to the beating relation} in \citet{Gonzalez-DiazHendrickxLohmann2014}.
Analogously, \citet{BoldiLuongoVigna2017} examine the rank monotonicity of centrality measures with respect to adding a new arc to a network.

A stronger version has been used in the axiomatic characterization of the ranking induced by the geometric mean \citep{Csato2018c}: $i \succeq j$ implies $i \succ j$ whenever $a_{ij}$ increases.
According to \citet[p.~188]{BouyssouPerny1992}, rank monotonicity seems an unobjectionable property in the context of partial ranking methods but its stronger version is much more demanding.

In some applications, the weights of the alternatives are more important than their ranking because, for example, certain resources will be allocated on the basis of the weights. Then it makes sense to exclude the possibility that the weight of any alternative decreases after one of its pairwise comparisons increases.

\begin{axiom} \label{Axiom2}
\emph{Weight monotonicity}:
Let $\mathbf{A} \in \mathcal{A}^{n \times n}$ be any pairwise comparison matrix and $1 \leq i,j \leq n$ be any two different alternatives. Let $\mathbf{A'} \in \mathcal{A}^{n \times n}$ be identical to $\mathbf{A}$ but $a_{ij}' > a_{ij}$ (and $a_{ji}' < a_{ji}$ due to the reciprocity property). 
The weighting method $f: \mathcal{A}^{n \times n} \to \mathcal{R}^n$ is called \emph{weight monotonic} if
\begin{equation} \label{eq_Landau}
a_{ij}' > a_{ij} \Rightarrow \frac{f_i \left( \mathbf{A'} \right)}{\sum_{k=1}^n f_k \left( \mathbf{A}' \right)} \geq \frac{f_i \left( \mathbf{A} \right)}{\sum_{k=1}^n f_k \left( \mathbf{A} \right)} \iff f_i \left( \mathbf{A'} \right) \geq f_i \left( \mathbf{A} \right),
\end{equation}
\end{axiom}

\citet[p.~201]{Landau1914} shows that weight monotonicity is violated by the principal right eigenvector for nonnegative matrices of order $n = 3$. The interesting and often neglected early observation is commented in \citet[p.~164]{DavidEdwards2001} with the following words: ``\emph{This leads Landau to dismiss the Kendall-Wei approach 40 years prior to its birth!}''
However, the set of pairwise comparison matrices is a strict subset of nonnegative matrices, and the example of Landau is not a pairwise comparison matrix.

Naturally, this type of monotonicity is a standard requirement in many fields of research such as cooperative game theory \citep{Young1985}, scientometrics \citep{KoczyStrobel2009, BouyssouMarchant2014, PerryReny2016}, the theoretical analysis of sports rules \citep{KendallLenten2017, VaziriDabadghaoYihMorin2018, Csato2019b, Csato2020f, Csato2020c}, or voting theory \citep{BalinskiYoung2001, Tasnadi2008}.

In the context of inconsistency indices, \citet{BrunelliFedrizzi2015} have suggested an axiom with the same flavour called \emph{monotonicity on single comparisons}, which formalizes a property proved by \citet{AupetitGenest1993} for the Saaty inconsistency index. The authors also provide a short overview of its origin in multi-criteria decision-making.

\subsection{Related works} \label{Sec23}

The current paper is not the first work focusing on a mathematical shortcoming of the eigenvector method.
\citet{JohnsonBeineWang1979} argue that the use of the left eigenvector is equally justified as long as the order is reversed, furthermore, the rankings from the two eigenvectors may disagree even when the matrix is nearly consistent. 
According to \citet{GenestLapointeDrury1993}, the ordering obtained from the principal right eigenvector depends on the choice of the parameter for numerically coded ordinal preferences.
\citet{BanaeCostaVansnick2008} find that the right eigenvector can violate a condition of order preservation, which is fundamental in decision aiding according to the authors. \citet{Kulakowski2015} examines the relationship between this property and the inconsistency index proposed by Saaty.
\citet{PerezMokotoff2016} present an example where the alternative with the highest priority for all decision-makers is not the best on the basis of their aggregated preferences. \citet{Csato2017b} traces back the origin of this problem to the right-left asymmetry \citep{JohnsonBeineWang1979}, and provides a minimal counterexample with four alternatives.

The eigenvector solution is not necessarily Pareto efficient, in other words, there may exist a weight vector which is at least as good in approximating all elements of the pairwise comparison matrix, and strictly better in at least one position \citep{BlanqueroCarrizosaConde2006}. \citet{Abele-NagyBozoki2016} prove that it is not possible if the pairwise comparison matrix differs from a consistent one only in one entry (and its reciprocal), while \citet{Abele-NagyBozokiRebak2018} extend this result to double perturbed matrices, which can be made consistent by altering two elements and their reciprocals. 
On the other hand, the eigenvector method may lead to an inefficient weight vector for matrices with an arbitrarily small inconsistency \citep{Bozoki2014}.
\citet{BozokiFulop2018} propose linear programs to test whether a given weight vector is efficient or not, and \citet{DulebaMoslem2019} give the first examination of this property on real data.

\section{Results} \label{Sec3}

The row geometric mean (logarithmic least squares) method trivially satisfies monotonicity: a greater value of $a_{ij}$ increases the weight of alternative $i$, decreases the weight of alternative $j$, while preserves the weights of all other alternatives before normalisation.
The case of the eigenvector method turns out to be more complicated.

\subsection{The geometric mean method is monotonic} \label{Sec31}

According to the following result, some weighting methods meet Axioms~\ref{Axiom1} and \ref{Axiom2}.

\begin{proposition} \label{Prop31}
The geometric mean and column sum methods satisfy both rank monotonicity and weight monotonicity.
\end{proposition}

\begin{proof}
Consider the pairwise comparison matrices $\mathbf{A}$ and $\mathbf{A'}$, which are identical except for $a_{ij}' > a_{ij}$ and $a_{ji}' < a_{ji}$.

\emph{Geometric mean method}: 
$i \succeq_{\mathbf{A}}^{GM} k$ means that $\prod_{\ell=1}^n a_{i \ell} \geq \prod_{\ell=1}^n a_{k \ell}$.
Since $\prod_{\ell=1}^n a_{i \ell}' > \prod_{\ell=1}^n a_{i \ell}$ and $\prod_{\ell=1}^n a_{k \ell}' \leq \prod_{\ell=1}^n a_{k \ell}$ if $k \neq i$, $\prod_{\ell=1}^n a_{i \ell}' > \prod_{\ell=1}^n a_{k \ell}'$, which implies $i \succ_{\mathbf{A'}}^{GM} k$, that is, rank monotonicity.
As only the weight of alternative $i$ increases, weight monotonicity holds, too.

\emph{Column sum method}: 
$i \succeq_{\mathbf{A}}^{CM} k$ means that $\sum_{\ell=1}^n a_{i \ell} \geq \sum_{\ell=1}^n a_{k \ell}$.
Since $\sum_{\ell=1}^n a_{i \ell}' > \sum_{\ell=1}^n a_{i \ell}$ and $\sum_{\ell=1}^n a_{k \ell}' \leq \sum_{\ell=1}^n a_{k \ell}$ if $k \neq i$, $\sum_{\ell=1}^n a_{i \ell}' > \sum_{\ell=1}^n a_{k \ell}'$, which implies $i \succ_{\mathbf{A'}}^{CM} k$, that is, rank monotonicity.
As only the weight of alternative $i$ increases, weight monotonicity holds, too.
\end{proof}

\citet[Proposition~2]{Csato2018c} has already shown the rank monotonicity of geometric mean.

\subsection{The eigenvector method violates monotonicity} \label{Sec32}

However, not all weighting methods are compatible with Axioms~\ref{Axiom1} and \ref{Axiom2}.

\begin{proposition} \label{Prop32}
The eigenvector method does not satisfy rank monotonicity.
\end{proposition}

\begin{proof}
It is sufficient to provide a counterexample.

\begin{example} \label{Examp31}
Consider the following parametric pairwise comparison matrix:
\[
\mathbf{A}^{\alpha} = \left[
\begin{array}{K{2.5em} K{2.5em} K{2.5em} K{2.5em} K{2.5em} K{2.5em}}
    1     & \alpha     & 8     &  1/9  &  1/4  &  1/9 \\
    1/ \alpha   & 1     &  1/5  &  1/9  & 1     &  1/4 \\
     1/8  & 5     & 1     &  1/8  &  1/2  &  1/7 \\
    9     & 9     & 8     & 1     & 7     & 8     \\
    4     & 1     & 2     &  1/7  & 1     &  1/9 \\
    9     & 4     & 7     &  1/8  & 9     & 1     \\
\end{array}
\right].
\]

\input{Figure1_EM_monotonicity_violation}

The weights of the first and the fifth alternatives are plotted in Figure~\ref{Fig1} as a function of parameter $\alpha$.
\end{example}

Note that the fifth alternative becomes more important than the first when $a_{12}$ is increased from $0.3$ to $0.5$, showing the violation of rank monotonicity by the eigenvector method.
\end{proof}

What does cause the violation of rank monotonicity in Example~\ref{Examp31}?
The growth of $\alpha$ increases $w_1$, which is also favourable for $w_5$ since $w_1 / w_5 \approx a_{15} = 1/4$. On the other hand, there is a ,,normalisation effect'', which decreases $w_5$ as the sum of all weights is fixed.

\begin{proposition} \label{Prop33}
The eigenvector method does not satisfy weight monotonicity.
\end{proposition}

\begin{proof}
Consider Example~\ref{Examp31}. According to Figure~\ref{Fig1}, the normalised weight of the first alternative has a minimum around $\alpha \approx 0.5$, thus the eigenvector method does not satisfy weight monotonicity.
\end{proof}

It remains difficult to explain the lack of weight monotonicity due to the multiple connections between the weights of alternatives as revealed by the formula $w_i^{EM} = \lambda_{\max} \sum_{\ell=1}^n a_{i \ell} w_{\ell}^{EM}$. Unfortunately, this is a common phenomenon for the unfavourable properties of the eigenvector method \citep{JohnsonBeineWang1979, BlanqueroCarrizosaConde2006, BanaeCostaVansnick2008}.

As it has already been mentioned, the principal right eigenvector violates Axiom~\ref{Axiom2} on the domain of nonnegative matrices of order $n = 3$ \citep[p.~201]{Landau1914}.
A section titled ``Non-monotonicity of some weight extraction methods'' on the \href{https://en.wikipedia.org/w/index.php?title=Analytic_hierarchy_process}{Wikipedia page of the Analytic Hierarchy Process} says that the eigenvector method is non-monotonic for reciprocal $n \times n$ matrices where $n > 3$. This contribution has been added on 10 December 2013 by the Dutch international draughts (dammen) player \href{https://toernooibase.kndb.nl/opvraag/liddetailp.php?SpId=2073&Id=f&taal=1}{\emph{Kees Pippel}}. We have managed to consult him, and he confirmed that it refers to the example of Landau, which does not imply Proposition~\ref{Prop33} due to the different domain.

\subsection{A framework for analysing the two axioms of monotonicity} \label{Sec33}

In order to investigate the probability of a non-monotonic eigenvector, we apply a computational technique and consider a large number of pairwise comparison matrices that are checked for both properties.

The entries of the random pairwise comparison matrices are generated according to the standard proposed by Saaty, which has also been used, for example, in \citet{BozokiRapcsak2008}.
Thus all entries $a_{ij}$ above the diagonal ($i < j$) are randomly chosen from the discrete set
\begin{equation} \label{eq_scale}
\left\{ \frac{1}{9};\, \frac{1}{8};\, \frac{1}{7};\, \frac{1}{6};\, \frac{1}{5};\, \frac{1}{4};\, \frac{1}{3};\, \frac{1}{2};\, 1;\, 2;\, 3;\, 4;\, 5;\, 6;\, 7;\, 8;\, 9 \right\}
\end{equation}
with equal probability $1/17$, and by setting $a_{ji} = 1 / a_{ij}$, as well as $a_{ii} = 1$.

\begin{remark} \label{Rem1}
The monotonicity of the principal right eigenvector is tested by substituting a matrix entry with the next element from the scale \eqref{eq_scale}, that is, if $2 \leq a \leq 9$, $a \in \mathbb{N}$, then $a_{ij} = 1/a$ is substituted with $1/(a-1)$, and if $1 \leq b \leq 9$, $b \in \mathbb{N}$, then $a_{ij} = b$ is substituted with $b+1$. This seems to be a realistic scenario since the decision-maker probably does not expect a counter-intuitive change in the ranking or weights when thinking over to give a more favourable opinion on an alternative.
\end{remark}

Consequently, one iteration of the computational process consists of the following steps:
\begin{enumerate}
\item \label{step1}
A random pairwise comparison matrix $\mathbf{A}$ of order $n$ is generated on the Saaty scale \eqref{eq_scale}.
\item
Its consistency ratio $CR(\mathbf{A})$ and eigenvector $\mathbf{w}^{EM}(\mathbf{A})$ is calculated, the number of matrices in the $m$th interval of consistency ratios, for which $0.01 (m-1) \leq CR(\mathbf{A}) < 0.01 m$, is increased by one.
\item
All entries above the diagonal are considered separately, hence $n(n-1)/2$ perturbed pairwise comparison matrices $\mathbf{A}^{ij}$ are introduced such that $\mathbf{A}^{ij} = \mathbf{A}$ except for its element in the $i$th row and $j$th column, $i < j$. This value is increased to $a_{ij}^{ij}$ following Remark~\ref{Rem1}, while reciprocity is preserved, thus $a_{ji}^{ij} = 1 / a_{ij}^{ij}$.
\item
The eigenvectors $\mathbf{w} \left( \mathbf{A}^{ij} \right)$ are computed.
\item
\emph{Monotonicity check}:
\begin{itemize}
\item
Rank monotonicity: $w^{EM}_i \left( \mathbf{A} \right) / w^{EM}_k \left( \mathbf{A} \right)$ and $w^{EM}_i \left( \mathbf{A}^{ij} \right) / w^{EM}_k \left( \mathbf{A}^{ij} \right)$ are compared for all $i < j$ and $1 \leq k \leq n$. \\
The $m$th interval of consistency ratios with the flag of rank monotonicity violation, in which $CR(\mathbf{A})$ falls, is increased by one if $w^{EM}_i(\mathbf{A}) > w^{EM}_k(\mathbf{A})$ but $w^{EM}_i \left( \mathbf{A}^{ij} \right) < w^{EM}_k \left( \mathbf{A}^{ij} \right)$, that is, Axiom~\ref{Axiom1} is not satisfied after increasing $a_{ij}$ on the scale \eqref{eq_scale}.
\item
Weight monotonicity: $w^{EM}_i \left( \mathbf{A} \right)$ and $w^{EM}_i \left( \mathbf{A}^{ij} \right)$ are compared for all $i < j$. \\
The $m$th interval of consistency ratios with the flag of weight monotonicity violation, in which $CR(\mathbf{A})$ falls, is increased by one if $w^{EM}_i(\mathbf{A}) > w^{EM}_i \left( \mathbf{A}^{ij} \right)$, that is, Axiom~\ref{Axiom2} is not satisfied after increasing $a_{ij}$ on the scale \eqref{eq_scale}.
\end{itemize}
\item \label{step6}
The pairwise comparison matrix $\mathbf{A}$, its consistency ratio $CR(\mathbf{A})$ and $i$, $j$, $k$ are saved as an example violating rank monotonicity, weight monotonicity, or both of them, if $CR(\mathbf{A})$ is smaller than the consistency ratio of all previously generated pairwise comparison matrices, where the eigenvector does not satisfy rank monotonicity, weight monotonicity, or both of them, respectively.
\end{enumerate}
Steps~\ref{step1}--\ref{step6} are repeated until the number of randomly generated pairwise comparison matrices reaches a predetermined limit.

\begin{table}[ht!]
\caption{Values of the random index used in the computations}
\centering
\label{Table1}
    \begin{tabularx}{\textwidth}{l CCCCCCC} \toprule
    Matrix size & 4     & 5     & 6     & 7     & 8     & 9     & 10 \\ \midrule
    Random index $RI_n$ & 0.884 & 1.109 & 1.249 & 1.341 & 1.404 & 1.451 & 1.486 \\ \bottomrule
    \end{tabularx}
\end{table}

The random indices $RI_n$ for the calculation of the consistency ratio $CR$ are reported in Table~\ref{Table1} for $4 \leq n \leq 10$. They are imported from \citet[Table~3]{BozokiRapcsak2008} and have been validated by our simulations, too. According to \citet{SaatyOzdemir2003}, it is reasonable to compare until nine items, hence we have not examined matrices with more than ten alternatives.

\subsection{The connection of monotonicity and inconsistency} \label{Sec34}

\begin{remark} \label{Rem2}
The eigenvector method satisfies both axioms of monotonicity if the number of alternatives is three because it is equivalent to the geometric mean method for $n = 3$ \citep{CrawfordWilliams1985}, and Proposition~\ref{Prop31} can be applied.
\end{remark}

For $n = 4$, there are only six elements above the diagonal, hence the total number of different matrices using the scale \eqref{eq_scale} is $17^6 = 24{,}137{,}569$. We have checked all of them through the process described in steps~\ref{step1}--\ref{step6}, and have not found any violation of monotonicity by the eigenvector method.
On the other hand, there are some cases where this problem emerges if the number of alternatives is at least five. Therefore, $10$ million randomly generated matrices have been examined for $5 \leq n \leq 8$, as well as $5$ million for $9 \leq n \leq 10$.

\input{Figure2_monotonicity_histogram}

Figure~\ref{Fig2} plots the proportion of pairwise comparison matrices for which the eigenvector method does not satisfy Axiom~\ref{Axiom1} or Axiom~\ref{Axiom2} as the function of the consistency ratio $CR$. Note that $CR$ cannot be arbitrarily large: \citet{AupetitGenest1993} derive a sharp upper bound on $\lambda_{\max}$ when the responses are coded on a bounded scale applied here.

The probability of violating rank monotonicity almost linearly increases when the pairwise comparison matrix becomes more inconsistent. There is no violation of Axiom~\ref{Axiom1} for nearly consistent matrices ($CI < 0.2$), and the problem emerges only with a probability of around 2\% even for heavily inconsistent matrices.

Surprisingly, this is not the case for weight monotonicity, Axiom~\ref{Axiom2} is violated primarily by matrices having a consistency ratio of about $0.5$. Weight monotonicity is satisfied more frequently than rank monotonicity for matrices of order $5$ and $6$, however, this result does not hold for moderately inconsistent matrices when $n \geq 7$.

\begin{table}[ht!]
\centering
\caption{The frequency of violating monotonicity properties by the eigenvector method}
\label{Table2}

\begin{subtable}{\textwidth}
\centering
\caption{All values of $CR$}
\label{Table2a}
\rowcolors{3}{gray!20}{}
    \begin{tabularx}{\textwidth}{cc CCc} \toprule \hiderowcolors
    Matrix size & Sample size & \multicolumn{2}{c}{Ratio of matrices} & Number of matrices violating  \\
          &       & Rank & Weight & both monotonicity properties \\ \bottomrule \showrowcolors   
    5     & $10^7$ & 0.47\% & 0.003\% & 0 \\
    6     & $10^7$ & 0.66\% & 0.070\% & 2 \\
    7     & $10^7$ & 0.80\% & 0.077\% & 11 \\
    8     & $10^7$ & 0.91\% & 0.046\% & 7 \\
    9     & $5 \times 10^6$ & 1.00\% & 0.020\% & 2 \\
    10    & $5 \times 10^6$ & 1.05\% & 0.007\% & 0 \\ \bottomrule
    \end{tabularx}
\end{subtable}

\vspace{0.5cm}
\begin{subtable}{\textwidth}
\centering
\caption{$CR < 0.5$}
\label{Table2b}
\rowcolors{3}{gray!20}{}
    \begin{tabularx}{\textwidth}{cr CCC} \toprule \hiderowcolors
    Matrix size & Sample size & \multicolumn{3}{c}{Number of matrices violating monotonicity} \\
          &       & Rank  & Weight & Both \\ \bottomrule \showrowcolors
    5     & 16,916,343 & 1,146  & 0     & 0 \\
    6     & 7,353,667 & 2,665  & 313   & 1 \\
    7     & 2,712,997 & 2,295  & 1,701  & 2 \\
    8     & 827,413 & 1,059  & 1,307  & 6 \\
    9     & 207,424 & 303   & 435   & 2 \\ 
    10    & 41,826 & 80    & 110   & 0 \\ \bottomrule
    \end{tabularx}
\end{subtable}
\end{table}

As Example~\ref{Examp31} shows, there are also some matrices for which the eigenvector method violates both rank and weight monotonicity simultaneously. Since these cases are rare, they are not plotted in Figure~\ref{Fig2} but summarised in Table~\ref{Table2}. The latter also reinforces the message of Figure~\ref{Fig2}: the violation of rank monotonicity becomes more serious as the number of alternatives increases, while weight monotonicity is the most problematic for $n = 7$ if one focuses on all randomly generated matrices. 

By limiting the monotonicity check to $CR < 0.5$, the number of randomly generated matrices can be increased to 100 million. However, only a fraction of them will be moderately inconsistent, see the column sample size in Table~\ref{Table2b}.
Then weight monotonicity becomes a more serious issue than rank monotonicity if $n \geq 8$. Furthermore, even the probability of violating weight monotonicity grows along with the size of the matrix for these moderately inconsistent matrices.

The two axioms are clearly independent as revealed by Figure~\ref{Fig2}. For large values of the consistency ratio $CR$, there is no violation of weight monotonicity (Axiom~\ref{Axiom2}) but there are several violations of rank monotonicity (Axiom~\ref{Axiom1}). On the other hand, it can be seen that weight monotonicity is violated more frequently compared to rank monotonicity for moderate levels of $CR$ if $n=9$. That is reinforced by Table~\ref{Table2b}, where there are more instances of violating rank monotonicity than weight monotonicity if $5 \leq n \leq 7$, however, this relation is reversed if $8 \leq n \leq 10$.
The relationship between the two properties is unsurprising in the view of their definitions. Firstly, it might happen that the (normalised) weight of alternative $i$ increases but the (normalised) weight of alternative $j$ grows faster after the entry $a_{ij}$ is increased. Secondly, it might happen that the (normalised) weight of alternative $i$ decreases if $a_{ij}$ is increased, without causing a rank reversal among the alternatives.

\section{Conclusions} \label{Sec4}

In the current paper, we have argued that two axioms of monotonicity on the numerical comparisons---called rank and weight monotonicity, respectively---are key requirements for any priority vector derived from a pairwise comparison matrix. The eigenvector method is proved to violate both properties.
However, contrary to the right-left asymmetry \citep{JohnsonBeineWang1979} and Pareto inefficiency \citep{Bozoki2014}, as well as, to the condition of order preservation \citep{BanaeCostaVansnick2008}, the emergence of a problematic situation seems to be avoidable for low levels of inconsistency.

Potential users should be informed about the rare occurrence of counter-intuitive changes in the weights after a matrix entry is increased. Nonetheless, it remains difficult to verify that there is no violation of rank or weight monotonicity in a specific application. Perhaps one can develop an online tool to check these properties, similarly to \href{http://pcmc.online/}{pcmc.online}, which is able to test the Pareto efficiency of the weights derived by the eigenvector method. The interaction with the decision-maker is indispensable: if the violation of monotonicity is judged problematic, the eigenvector method should be exchanged for the row geometric mean. It seems also probable that the violation of rank or weight monotonicity indicates some discrepancies in the data as shown by the lack of the problem at low levels of inconsistency. For instance, the fourth alternative dominates in Example~\ref{Examp31}.

To summarise, our findings have important implications for practitioners.
First, they emphasise the need for inconsistency reduction methods \citep{AbelMikhailovKeane2018, AguaronEscobarMoreno-Jimenez2020, BozokiFulopPoesz2015, ErguKouPengShi2011}.
Second, the possibly strange behaviour of the right eigenvector makes the use of this method questionable for inherently inconsistent matrices such as the ones emerging in sports applications \citep{BozokiCsatoTemesi2016, ChaoKouLiPeng2018}, where rewarding players or teams for poor performance is unfair \citep{KendallLenten2017, VaziriDabadghaoYihMorin2018}.
Third, the violation of monotonicity by the eigenvector method serves as a powerful argument for preferring the use of the geometric mean (logarithmic least squares) method in applications of AHP.

\section*{Acknowledgements}
\addcontentsline{toc}{section}{Acknowledgements}
\noindent
\emph{L\'aszl\'o Csat\'o}, the father of the first author has made a substantial contribution to the paper by helping to code the computations in Python. \\
We are grateful to \emph{S\'andor Boz\'oki}, \emph{J\'anos F\"ul\"op}, \emph{Tam\'as Halm}, \emph{Kees Pippel}, and \emph{Bal\'azs Vida} for useful advice. \\
The research was supported by the MTA Premium Postdoctoral Research Program grant PPD2019-9/2019, the NKFIH grant K 128573, and the \'UNKP-19-3-III-BCE-97 New National Excellence Program of the Ministry for Innovation and Technology.

\bibliographystyle{apalike}
\bibliography{All_references}

\end{document}

%% file: Figure1_EM_monotonicity_violation.tex
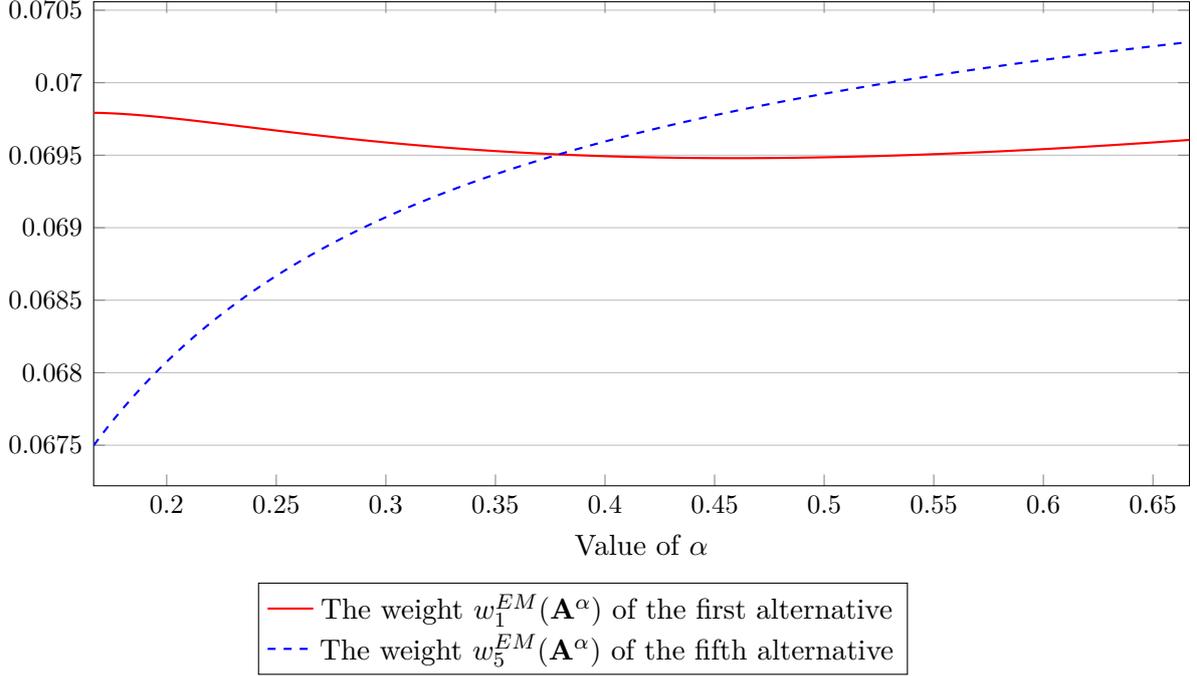
\begin{figure}[t]
\centering
\caption{Weights derived by the eigenvector method in Example~\ref{Examp31}}
\label{Fig1}

\begin{tikzpicture}
\begin{axis}[
width = 1\textwidth,
height = 0.5\textwidth,
ymajorgrids = true,
xmin = 0.1666666665,
xmax = 0.666666666,
xlabel = {Value of $\alpha$},
xlabel style = {font=\small},
scaled ticks = false,
tick label style = {/pgf/number format/.cd,fixed,precision=4},
legend style = {font=\small,at={(0.15,-0.2)},anchor=north west,legend columns=1},
legend entries = {The weight $w_1^{EM}(\mathbf{A}^{\alpha})$ of the first alternative,The weight $w_5^{EM}(\mathbf{A}^{\alpha})$ of the fifth alternative},
]      

\addplot [red, thick] coordinates {
(0.1666666665,0.0697925958103146)
(0.16999999983,0.0697915174877109)
(0.17333333316,0.0697897636916249)
(0.17666666649,0.0697874065650183)
(0.17999999982,0.0697845116032374)
(0.18333333315,0.0697811383099114)
(0.18666666648,0.0697773407813042)
(0.18999999981,0.0697731682277667)
(0.19333333314,0.0697686654397755)
(0.19666666647,0.0697638732050768)
(0.1999999998,0.0697588286825962)
(0.20333333313,0.0697535657380652)
(0.20666666646,0.0697481152456857)
(0.20999999979,0.0697425053596208)
(0.21333333312,0.0697367617586327)
(0.21666666645,0.0697309078667934)
(0.21999999978,0.069724965052836)
(0.22333333311,0.0697189528104185)
(0.22666666644,0.0697128889213051)
(0.22999999977,0.0697067896032379)
(0.2333333331,0.0697006696440714)
(0.23666666643,0.0696945425235672)
(0.23999999976,0.0696884205240869)
(0.24333333309,0.0696823148312892)
(0.24666666642,0.0696762356258142)
(0.24999999975,0.0696701921668343)
(0.25333333308,0.0696641928682547)
(0.25666666641,0.0696582453682683)
(0.25999999974,0.0696523565928914)
(0.26333333307,0.0696465328140473)
(0.2666666664,0.0696407797027009)
(0.26999999973,0.0696351023775038)
(0.27333333306,0.0696295054493553)
(0.27666666639,0.0696239930622528)
(0.27999999972,0.069618568930761)
(0.28333333305,0.0696132363744023)
(0.28666666638,0.0696079983492394)
(0.28999999971,0.0696028574768965)
(0.29333333304,0.0695978160712383)
(0.29666666637,0.0695928761629142)
(0.2999999997,0.0695880395219436)
(0.30333333303,0.0695833076785125)
(0.30666666636,0.0695786819421317)
(0.30999999969,0.0695741634192911)
(0.31333333302,0.0695697530297373)
(0.31666666635,0.0695654515214875)
(0.31999999968,0.0695612594846812)
(0.32333333301,0.0695571773643675)
(0.32666666634,0.0695532054723096)
(0.32999999967,0.0695493439978912)
(0.333333333,0.0695455930181887)
(0.33666666633,0.0695419525072844)
(0.33999999966,0.0695384223448716)
(0.34333333299,0.0695350023242123)
(0.34666666632,0.0695316921595)
(0.34999999965,0.0695284914926642)
(0.35333333298,0.0695253998996715)
(0.35666666631,0.0695224168963527)
(0.35999999964,0.0695195419437975)
(0.36333333297,0.069516774453346)
(0.3666666663,0.0695141137912099)
(0.36999999963,0.0695115592827496)
(0.37333333296,0.0695091102164345)
(0.37666666629,0.0695067658475079)
(0.37999999962,0.0695045254013812)
(0.38333333295,0.0695023880767742)
(0.38666666628,0.0695003530486229)
(0.38999999961,0.0694984194707697)
(0.39333333294,0.0694965864784546)
(0.39666666627,0.0694948531906172)
(0.3999999996,0.0694932187120296)
(0.40333333293,0.0694916821352658)
(0.40666666626,0.0694902425425241)
(0.40999999959,0.0694888990073122)
(0.41333333292,0.0694876505960009)
(0.41666666625,0.0694864963692628)
(0.41999999958,0.0694854353833982)
(0.42333333291,0.0694844666915572)
(0.42666666624,0.0694835893448701)
(0.42999999957,0.0694828023934849)
(0.4333333329,0.0694821048875256)
(0.43666666623,0.0694814958779714)
(0.43999999956,0.0694809744174677)
(0.44333333289,0.0694805395610679)
(0.44666666622,0.0694801903669151)
(0.44999999955,0.0694799258968663)
(0.45333333288,0.0694797452170644)
(0.45666666621,0.0694796473984593)
(0.45999999954,0.0694796315172832)
(0.46333333287,0.0694796966554846)
(0.4666666662,0.0694798419011209)
(0.46999999953,0.0694800663487145)
(0.47333333286,0.0694803690995756)
(0.47666666619,0.0694807492620917)
(0.47999999952,0.0694812059519879)
(0.48333333285,0.0694817382925607)
(0.48666666618,0.0694823454148838)
(0.48999999951,0.069483026457993)
(0.49333333284,0.069483780569047)
(0.49666666617,0.0694846069034684)
(0.4999999995,0.0694855046250663)
(0.50333333283,0.0694864729061397)
(0.50666666616,0.069487510927566)
(0.50999999949,0.0694886178788745)
(0.51333333282,0.069489792958304)
(0.51666666615,0.0694910353728496)
(0.51999999948,0.0694923443382962)
(0.52333333281,0.0694937190792419)
(0.52666666614,0.0694951588291089)
(0.52999999947,0.0694966628301484)
(0.5333333328,0.0694982303334336)
(0.53666666613,0.0694998605988464)
(0.53999999946,0.069501552895056)
(0.54333333279,0.0695033064994912)
(0.54666666612,0.0695051206983059)
(0.54999999945,0.0695069947863385)
(0.55333333278,0.0695089280670679)
(0.55666666611,0.0695109198525623)
(0.55999999944,0.0695129694634249)
(0.56333333277,0.069515076228736)
(0.5666666661,0.0695172394859904)
(0.56999999943,0.0695194585810315)
(0.57333333276,0.0695217328679842)
(0.57666666609,0.0695240617091832)
(0.57999999942,0.069526444475099)
(0.58333333275,0.0695288805442631)
(0.58666666608,0.0695313693031903)
(0.58999999941,0.0695339101462999)
(0.59333333274,0.0695365024758339)
(0.59666666607,0.0695391457017763)
(0.5999999994,0.0695418392417701)
(0.60333333273,0.0695445825210318)
(0.60666666606,0.0695473749722691)
(0.60999999939,0.0695502160355932)
(0.61333333272,0.0695531051584347)
(0.61666666605,0.0695560417954561)
(0.61999999938,0.0695590254084665)
(0.62333333271,0.0695620554663337)
(0.62666666604,0.0695651314448984)
(0.62999999937,0.0695682528268872)
(0.6333333327,0.0695714191018257)
(0.63666666603,0.0695746297659531)
(0.63999999936,0.0695778843221346)
(0.64333333269,0.0695811822797772)
(0.64666666602,0.0695845231547428)
(0.64999999935,0.0695879064692648)
(0.65333333268,0.069591331751862)
(0.65666666601,0.0695947985372558)
(0.65999999934,0.0695983063662868)
(0.66333333267,0.0696018547858315)
(0.666666666,0.0696054433487204)
};

\addplot [blue, thick, dashed] coordinates {
(0.1666666665,0.0674983769809361)
(0.16999999983,0.0675655400366738)
(0.17333333316,0.0676303113583927)
(0.17666666649,0.0676928155422652)
(0.17999999982,0.0677531688218542)
(0.18333333315,0.067811479752052)
(0.18666666648,0.0678678498272221)
(0.18999999981,0.0679223740408118)
(0.19333333314,0.0679751413927918)
(0.19666666647,0.0680262353505072)
(0.1999999998,0.0680757342678476)
(0.20333333313,0.0681237117670613)
(0.20666666646,0.0681702370870356)
(0.20999999979,0.0682153754014223)
(0.21333333312,0.0682591881096025)
(0.21666666645,0.068301733103157)
(0.21999999978,0.0683430650102019)
(0.22333333311,0.0683832354197047)
(0.22666666644,0.0684222930876616)
(0.22999999977,0.0684602841268181)
(0.2333333331,0.0684972521814415)
(0.23666666643,0.0685332385884982)
(0.23999999976,0.0685682825264458)
(0.24333333309,0.0686024211527376)
(0.24666666642,0.0686356897310165)
(0.24999999975,0.0686681217488908)
(0.25333333308,0.0686997490270893)
(0.25666666641,0.0687306018207196)
(0.25999999974,0.0687607089132889)
(0.26333333307,0.0687900977040778)
(0.2666666664,0.0688187942894057)
(0.26999999973,0.0688468235382835)
(0.27333333306,0.0688742091628917)
(0.27666666639,0.0689009737842959)
(0.27999999972,0.0689271389937653)
(0.28333333305,0.0689527254100347)
(0.28666666638,0.0689777527328164)
(0.28999999971,0.0690022397928451)
(0.29333333304,0.0690262045987128)
(0.29666666637,0.0690496643807305)
(0.2999999997,0.0690726356320349)
(0.30333333303,0.0690951341471356)
(0.30666666636,0.0691171750580894)
(0.30999999969,0.0691387728684665)
(0.31333333302,0.0691599414852666)
(0.31666666635,0.0691806942489224)
(0.31999999968,0.0692010439615279)
(0.32333333301,0.0692210029134076)
(0.32666666634,0.0692405829081402)
(0.32999999967,0.0692597952861411)
(0.333333333,0.0692786509468967)
(0.33666666633,0.0692971603699424)
(0.33999999966,0.0693153336346622)
(0.34333333299,0.0693331804389877)
(0.34666666632,0.0693507101170694)
(0.34999999965,0.0693679316559786)
(0.35333333298,0.0693848537115086)
(0.35666666631,0.0694014846231254)
(0.35999999964,0.0694178324281236)
(0.36333333297,0.0694339048750347)
(0.3666666663,0.0694497094363332)
(0.36999999963,0.0694652533204846)
(0.37333333296,0.0694805434833704)
(0.37666666629,0.0694955866391331)
(0.37999999962,0.0695103892704686)
(0.38333333295,0.0695249576384038)
(0.38666666628,0.0695392977915857)
(0.38999999961,0.0695534155751104)
(0.39333333294,0.0695673166389209)
(0.39666666627,0.0695810064457938)
(0.3999999996,0.0695944902789419)
(0.40333333293,0.06960777324925)
(0.40666666626,0.0696208603021695)
(0.40999999959,0.0696337562242848)
(0.41333333292,0.0696464656495727)
(0.41666666625,0.0696589930653706)
(0.41999999958,0.0696713428180696)
(0.42333333291,0.0696835191185458)
(0.42666666624,0.0696955260473443)
(0.42999999957,0.0697073675596293)
(0.4333333329,0.0697190474899129)
(0.43666666623,0.0697305695565713)
(0.43999999956,0.0697419373661626)
(0.44333333289,0.0697531544175554)
(0.44666666622,0.0697642241058754)
(0.44999999955,0.0697751497262831)
(0.45333333288,0.0697859344775874)
(0.45666666621,0.0697965814657058)
(0.45999999954,0.0698070937069773)
(0.46333333287,0.0698174741313354)
(0.4666666662,0.0698277255853492)
(0.46999999953,0.0698378508351366)
(0.47333333286,0.0698478525691581)
(0.47666666619,0.0698577334008955)
(0.47999999952,0.0698674958714219)
(0.48333333285,0.0698771424518663)
(0.48666666618,0.0698866755457805)
(0.48999999951,0.0698960974914103)
(0.49333333284,0.0699054105638766)
(0.49666666617,0.0699146169772702)
(0.4999999995,0.0699237188866651)
(0.50333333283,0.0699327183900521)
(0.50666666616,0.0699416175301985)
(0.50999999949,0.0699504182964358)
(0.51333333282,0.0699591226263779)
(0.51666666615,0.0699677324075756)
(0.51999999948,0.0699762494791076)
(0.52333333281,0.069984675633112)
(0.52666666614,0.0699930126162589)
(0.52999999947,0.0700012621311715)
(0.5333333328,0.0700094258377917)
(0.53666666613,0.070017505354697)
(0.53999999946,0.0700255022603693)
(0.54333333279,0.0700334180944166)
(0.54666666612,0.070041254358752)
(0.54999999945,0.0700490125187296)
(0.55333333278,0.0700566940042397)
(0.55666666611,0.0700643002107658)
(0.55999999944,0.0700718325004036)
(0.56333333277,0.0700792922028449)
(0.5666666661,0.0700866806163266)
(0.56999999943,0.0700939990085471)
(0.57333333276,0.070101248617551)
(0.57666666609,0.0701084306525838)
(0.57999999942,0.0701155462949162)
(0.58333333275,0.0701225966986423)
(0.58666666608,0.070129582991449)
(0.58999999941,0.0701365062753614)
(0.59333333274,0.0701433676274609)
(0.59666666607,0.0701501681005824)
(0.5999999994,0.0701569087239855)
(0.60333333273,0.0701635905040061)
(0.60666666606,0.0701702144246856)
(0.60999999939,0.0701767814483784)
(0.61333333272,0.0701832925163428)
(0.61666666605,0.0701897485493096)
(0.61999999938,0.0701961504480348)
(0.62333333271,0.0702024990938337)
(0.62666666604,0.0702087953490981)
(0.62999999937,0.0702150400577972)
(0.6333333327,0.070221234045963)
(0.63666666603,0.0702273781221608)
(0.63999999936,0.0702334730779438)
(0.64333333269,0.0702395196882954)
(0.64666666602,0.0702455187120565)
(0.64999999935,0.0702514708923407)
(0.65333333268,0.0702573769569358)
(0.65666666601,0.0702632376186939)
(0.65999999934,0.0702690535759105)
(0.66333333267,0.0702748255126893)
(0.666666666,0.0702805540993)
};
\end{axis}
\end{tikzpicture}

\end{figure}


%% file: Figure2_monotonicity_histogram.tex
\begin{figure}[ht!]
\centering
\caption{Randomly generated pairwise comparison matrices for which \\ the eigenvector method does not satisfy monotonicity}
\label{Fig2}

\begin{tikzpicture}
\begin{axis}[
name = axis1,
title = {Matrix size: $n=5$; sample size: $10^7$},
title style = {font = \small},
xlabel = Value of $CR$,
xlabel style = {font = \small},
width = 0.48\textwidth,
height = 0.35\textwidth,
ymajorgrids = true,
xmin = 0,
xmax = 2.4,
ymin = -0.001,
ymax = 0.021,
scaled y ticks = false,
y tick label style = {/pgf/number format/.cd,fixed,precision=4},
]      

\addplot [black, dotted, very thick] coordinates {
(0,0)
(0.1,0)
(0.2,0)
(0.3,2.29977136439686E-05)
(0.4,0.000167819750760602)
(0.5,0.000587831198532161)
(0.6,0.00121173728366185)
(0.7,0.00221660708531227)
(0.8,0.00275743552879863)
(0.9,0.00423778599047232)
(1,0.00467502805931458)
(1.1,0.00676029717854602)
(1.2,0.00713598037882845)
(1.3,0.0106409357330076)
(1.4,0.00911355234524527)
(1.5,0.00811091956812096)
(1.6,0.00793885748895268)
(1.7,0.00759498381274886)
(1.8,0.00817402476574441)
(1.9,0.00891835859307978)
(2,0.00887470071827614)
(2.1,0.008934575847824)
(2.2,0.010789567702751)
(2.3,0.0112014453477868)
(2.4,0.00988286969253294)
};

\addplot [blue, dashdotted, very thick] coordinates {
(0,0)
(0.1,0)
(0.2,0)
(0.3,0)
(0.4,0)
(0.5,1.73914555778746E-06)
(0.6,5.03536940931868E-05)
(0.7,0.000188230900586979)
(0.8,0.000150226884715608)
(0.9,6.54226570677832E-05)
(1,3.78915074678935E-05)
(1.1,3.95956492300626E-06)
(1.2,0)
(1.3,0)
(1.4,0)
(1.5,0)
(1.6,0)
(1.7,0)
(1.8,0)
(1.9,0)
(2,0)
(2.1,0)
(2.2,0)
(2.3,0)
(2.4,0)
};
\legend{}
\end{axis}

\begin{axis}[
width = 0.48\textwidth, 
height = 0.35\textwidth,
xmin = 0,
xmax = 2.4,
xticklabels={,,},
axis y line* = right,
ymode = log,
ymin = 1,
ymax = 5000000,
]
\addplot[red,loosely dashed, very thick] coordinates {
(0,23853)
(0.1,187304)
(0.2,404135)
(0.3,521791)
(0.4,554166)
(0.5,574995)
(0.6,615645)
(0.7,664078)
(0.8,712256)
(0.9,748976)
(1,765343)
(1.1,757659)
(1.2,720854)
(1.3,657367)
(1.4,572115)
(1.5,473929)
(1.6,368189)
(1.7,268730)
(1.8,180572)
(1.9,111904)
(2,62650)
(2.1,31227)
(2.2,13995)
(2.3,5535)
(2.4,2732)
};
\end{axis}

\begin{axis}[
at = {(axis1.south east)},
xshift = 0.12\textwidth,
title = {Matrix size: $n=6$; sample size: $10^7$},
title style = {font = \small},
xlabel = Value of $CR$,
xlabel style = {font = \small},
width = 0.48\textwidth,
height = 0.35\textwidth,
ymajorgrids = true,
xmin = 0,
xmax = 2,
ymin = -0.001,
ymax = 0.021,
scaled y ticks = false,
y tick label style = {/pgf/number format/.cd,fixed,precision=4},
]      

\addplot [black, dotted, very thick] coordinates {
(0,0)
(0.1,0)
(0.2,0)
(0.3,0.000157226023031488)
(0.4,0.000696977246072155)
(0.5,0.00158218305073842)
(0.6,0.00268950055146062)
(0.7,0.00408486530647193)
(0.8,0.00578378753757478)
(0.9,0.00689540885051115)
(1,0.00722296076626785)
(1.1,0.00783917090319332)
(1.2,0.00866433041816108)
(1.3,0.00908975488212282)
(1.4,0.00993711295665115)
(1.5,0.0110368557295538)
(1.6,0.0114543711537867)
(1.7,0.0124833625218914)
(1.8,0.0122887284001966)
(1.9,0.0144272989587254)
(2,0.0114660114660115)
};

\addplot [blue, dashdotted, very thick] coordinates {
(0,0)
(0.1,0)
(0.2,0)
(0.3,0)
(0.4,0.000131557954077877)
(0.5,0.000842356668690467)
(0.6,0.00153384401198146)
(0.7,0.00166962180212128)
(0.8,0.00133632278097997)
(0.9,0.000873778836473974)
(1,0.000517264887343628)
(1.1,0.000324367194357437)
(1.2,0.000136680814329904)
(1.3,8.06743849570045E-05)
(1.4,3.45638711535692E-05)
(1.5,1.28149268267678E-05)
(1.6,0)
(1.7,0)
(1.8,0)
(1.9,0)
(2,0)
};
\legend{}
\end{axis}

\begin{axis}[
at = {(axis1.south east)},
xshift = 0.12\textwidth,
width = 0.48\textwidth, 
height = 0.35\textwidth,
xmin = 0,
xmax = 2,
xticklabels={,,},
axis y line* = right,
ymode = log,
ymin = 1,
ymax = 5000000,
]
\addplot[red,loosely dashed, very thick] coordinates {
(0,811)
(0.1,25252)
(0.2,115221)
(0.3,235330)
(0.4,357257)
(0.5,497414)
(0.6,663692)
(0.7,840909)
(0.8,1007990)
(0.9,1128432)
(1,1173480)
(1.1,1122185)
(1.2,973070)
(1.3,756126)
(1.4,520775)
(1.5,312136)
(1.6,161685)
(1.7,71375)
(1.8,26447)
(1.9,7971)
(2,2442)
};
\end{axis}
\end{tikzpicture}

\vspace{0.5cm}
\begin{tikzpicture}
\begin{axis}[
name = axis3,
title = {Matrix size: $n=7$; sample size: $10^7$},
title style = {font = \small},
xlabel = Value of $CR$,
xlabel style = {font = \small},
width = 0.48\textwidth,
height = 0.35\textwidth,
ymajorgrids = true,
xmin = 0,
xmax = 1.8,
ymin = -0.001,
ymax = 0.021,
scaled y ticks = false,
y tick label style = {/pgf/number format/.cd,fixed,precision=4},
]      

\addplot [black, dotted, very thick] coordinates {
(0,0)
(0.1,0)
(0.2,4.85107208693121E-05)
(0.3,0.000307815845824411)
(0.4,0.00121075746586178)
(0.5,0.00237801629882953)
(0.6,0.00388784306631008)
(0.7,0.00507664676214651)
(0.8,0.00613449623953854)
(0.9,0.00722249036028916)
(1,0.00839000605379664)
(1.1,0.00946888491728257)
(1.2,0.0107101126947489)
(1.3,0.0117096092641771)
(1.4,0.0128467267177458)
(1.5,0.0134589732188665)
(1.6,0.0143432528758946)
(1.7,0.013998133582189)
(1.8,0.0180932498260264)
};

\addplot [blue, dashdotted, very thick] coordinates {
(0,0)
(0.1,0)
(0.2,0)
(0.3,9.36830835117773E-05)
(0.4,0.000936623700006282)
(0.5,0.00235152849352079)
(0.6,0.00264025163458371)
(0.7,0.00187778492072434)
(0.8,0.00106672921384145)
(0.9,0.000580106394827705)
(1,0.000255295808219767)
(1.1,0.000116742584723278)
(1.2,3.92031622609324E-05)
(1.3,4.73370540257798E-06)
(1.4,3.25068995894379E-06)
(1.5,0)
(1.6,0)
(1.7,0)
(1.8,0)
};
\legend{}
\end{axis}

\begin{axis}[
width = 0.48\textwidth, 
height = 0.35\textwidth,
xmin = 0,
xmax = 1.8,
xticklabels={,,},
axis y line* = right,
ymode = log,
ymin = 1,
ymax = 5000000,
]
\addplot[red,loosely dashed, very thick] coordinates {
(0,10)
(0.1,1803)
(0.2,20614)
(0.3,74720)
(0.4,175097)
(0.5,339779)
(0.6,585929)
(0.7,905322)
(0.8,1244927)
(0.9,1508344)
(1,1582478)
(1.1,1413366)
(1.2,1045834)
(1.3,633753)
(1.4,307627)
(1.5,116948)
(1.6,34511)
(1.7,7501)
(1.8,1437)
};
\end{axis}

\begin{axis}[
at = {(axis3.south east)},
xshift = 0.12\textwidth,
title = {Matrix size: $n=8$; sample size: $10^7$},
title style = {font = \small},
xlabel = Value of $CR$,
xlabel style = {font = \small},
width = 0.48\textwidth,
height = 0.35\textwidth,
ymajorgrids = true,
xmin = 0,
xmax = 1.6,
ymin = -0.001,
ymax = 0.021,
scaled y ticks = false,
y tick label style = {/pgf/number format/.cd,fixed,precision=4},
max space between ticks = 40,
]      

\addplot [black, dotted, very thick] coordinates {
(0.1,0)
(0.2,0)
(0.3,0.00036138047340842)
(0.4,0.00159123882622736)
(0.5,0.00271493702103265)
(0.6,0.0039922021259898)
(0.7,0.00521626838139747)
(0.8,0.00655626886642026)
(0.9,0.00804711064604157)
(1,0.00949559885323552)
(1.1,0.0109304906421405)
(1.2,0.0123597701113775)
(1.3,0.0140709758297892)
(1.4,0.0147744945567652)
(1.5,0.0160297203398576)
(1.6,0.0196581196581197)
};

\addplot [blue, dashdotted, very thick] coordinates {
(0.1,0)
(0.2,0)
(0.3,0.00018069023670421)
(0.4,0.00184084491661597)
(0.5,0.00298534907492036)
(0.6,0.00227248428710189)
(0.7,0.00119549109272506)
(0.8,0.000601038453026016)
(0.9,0.00025494847539215)
(1,8.20419740919549E-05)
(1.1,2.62887594897727E-05)
(1.2,7.23700466580029E-06)
(1.3,0)
(1.4,0)
(1.5,0)
(1.6,0)
};
\legend{}
\end{axis}

\begin{axis}[
at = {(axis3.south east)},
xshift = 0.12\textwidth,
width = 0.48\textwidth, 
height = 0.35\textwidth,
xmin = 0,
xmax = 1.6,
xticklabels={,,},
axis y line* = right,
ymode = log,
ymin = 1,
ymax = 5000000,
]
\addplot[red,loosely dashed, very thick] coordinates {
(0.1,71)
(0.2,2181)
(0.3,16603)
(0.4,64101)
(0.5,184903)
(0.6,439607)
(0.7,865753)
(0.8,1414219)
(0.9,1878811)
(1,1974599)
(1.1,1597641)
(1.2,967251)
(1.3,426765)
(1.4,133744)
(1.5,29071)
(1.6,4680)
};
\end{axis}
\end{tikzpicture}

\vspace{0.5cm}
\begin{tikzpicture}
\begin{axis}[
name = axis5,
title = {Matrix size: $n=9$; sample size: $5 \times 10^6$},
title style = {font = \small},
xlabel = Value of $CR$,
xlabel style = {font = \small},
width = 0.48\textwidth,
height = 0.35\textwidth,
ymajorgrids = true,
xmin = 0,
xmax = 1.5,
ymin = -0.001,
ymax = 0.021,
scaled y ticks = false,
y tick label style = {/pgf/number format/.cd,fixed,precision=4},
]      

\addplot [black, dotted, very thick] coordinates {
(0.2,0)
(0.3,0)
(0.4,0.00118152524167562)
(0.5,0.00305560233861848)
(0.6,0.00421433713299274)
(0.7,0.00577505906310405)
(0.8,0.00721926793163418)
(0.9,0.00861465241253144)
(1,0.0103735644452034)
(1.1,0.0119306491462499)
(1.2,0.0137206269271073)
(1.3,0.015878741108576)
(1.4,0.0170363083527097)
(1.5,0.0156608097784568)
};

\addplot [blue, dashdotted, very thick] coordinates {
(0.2,0)
(0.3,0.000742942050520059)
(0.4,0.00225563909774436)
(0.5,0.0025503452590044)
(0.6,0.0013908011433923)
(0.7,0.000734477208668422)
(0.8,0.000274667592789447)
(0.9,8.03938943517036E-05)
(1,2.81774324382018E-05)
(1.1,4.77655856120503E-06)
(1.2,0)
(1.3,0)
(1.4,0)
(1.5,0)
};
\legend{}
\end{axis}

\begin{axis}[
width = 0.48\textwidth, 
height = 0.35\textwidth,
xmin = 0,
xmax = 1.5,
xticklabels={,,},
axis y line* = right,
ymode = log,
ymin = 1,
ymax = 5000000,
]
\addplot[red,loosely dashed, very thick] coordinates {
(0.2,80)
(0.3,1346)
(0.4,9310)
(0.5,41563)
(0.6,143083)
(0.7,377139)
(0.8,757279)
(0.9,1119488)
(1,1171150)
(1.1,837423)
(1.2,397941)
(1.3,119216)
(1.4,22364)
(1.5,2618)
};
\end{axis}

\begin{axis}[
at = {(axis5.south east)},
xshift = 0.12\textwidth,
title = {Matrix size: $n=10$; sample size: $5 \times 10^6$},
title style = {font = \small},
xlabel = Value of $CR$,
xlabel style = {font = \small},
width = 0.48\textwidth,
height = 0.35\textwidth,
ymajorgrids = true,
xmin = 0,
xmax = 1.4,
ymin = -0.001,
ymax = 0.021,
scaled y ticks = false,
y tick label style = {/pgf/number format/.cd,fixed,precision=4},
max space between ticks = 40,
]      

\addplot [black, dotted, very thick] coordinates {
(0.2,0)
(0.3,0)
(0.4,0.000527148128624143)
(0.5,0.00287769784172662)
(0.6,0.00423250564334086)
(0.7,0.00582254136416659)
(0.8,0.0073779152406128)
(0.9,0.00912766187775942)
(1,0.0110339528577715)
(1.1,0.0127072416715406)
(1.2,0.0145288901136481)
(1.3,0.0162591896050607)
(1.4,0.0191372549019608)
};

\addplot [blue, dashdotted, very thick] coordinates {
(0.2,0)
(0.3,0)
(0.4,0.000527148128624143)
(0.5,0.00196206671026815)
(0.6,0.00101825498086171)
(0.7,0.000317232673643582)
(0.8,9.25007914680397E-05)
(0.9,2.57007544339642E-05)
(1,4.46718739181031E-06)
(1.1,2.38723307750152E-06)
(1.2,0)
(1.3,0)
(1.4,0)
};
\legend{}
\end{axis}

\begin{axis}[
at = {(axis3.south east)},
xshift = 0.12\textwidth,
width = 0.48\textwidth, 
height = 0.35\textwidth,
xmin = 0,
xmax = 1.4,
xticklabels={,,},
axis y line* = right,
ymode = log,
ymin = 1,
ymax = 5000000,
]
\addplot[red,loosely dashed, very thick] coordinates {
(0.2,1)
(0.3,138)
(0.4,1897)
(0.5,15290)
(0.6,81512)
(0.7,302617)
(0.8,767561)
(0.9,1284009)
(1,1343127)
(1.1,837790)
(1.2,301193)
(1.3,58490)
(1.4,6375)
};
\end{axis}
\end{tikzpicture}

\vspace{-0.2cm}
\begin{center}
\begin{tikzpicture}
	\begin{customlegend}[legend columns=1,legend entries={Ratio of matrices exhibiting a violation of rank monotonicity (left scale)$\quad$, Ratio of matrices exhibiting a violation of weight monotonicity (left scale), Number of all randomly generated matrices (right log scale)$\qquad \qquad \qquad$}, legend style={font = \small}]
        \addlegendimage{color = black, dotted, very thick}
        \addlegendimage{color = blue, dashdotted, very thick}
        \addlegendimage{color = red, loosely dashed, very thick} 
	\end{customlegend}
\end{tikzpicture}
\end{center}
\vspace{-0.5cm}
\end{figure}
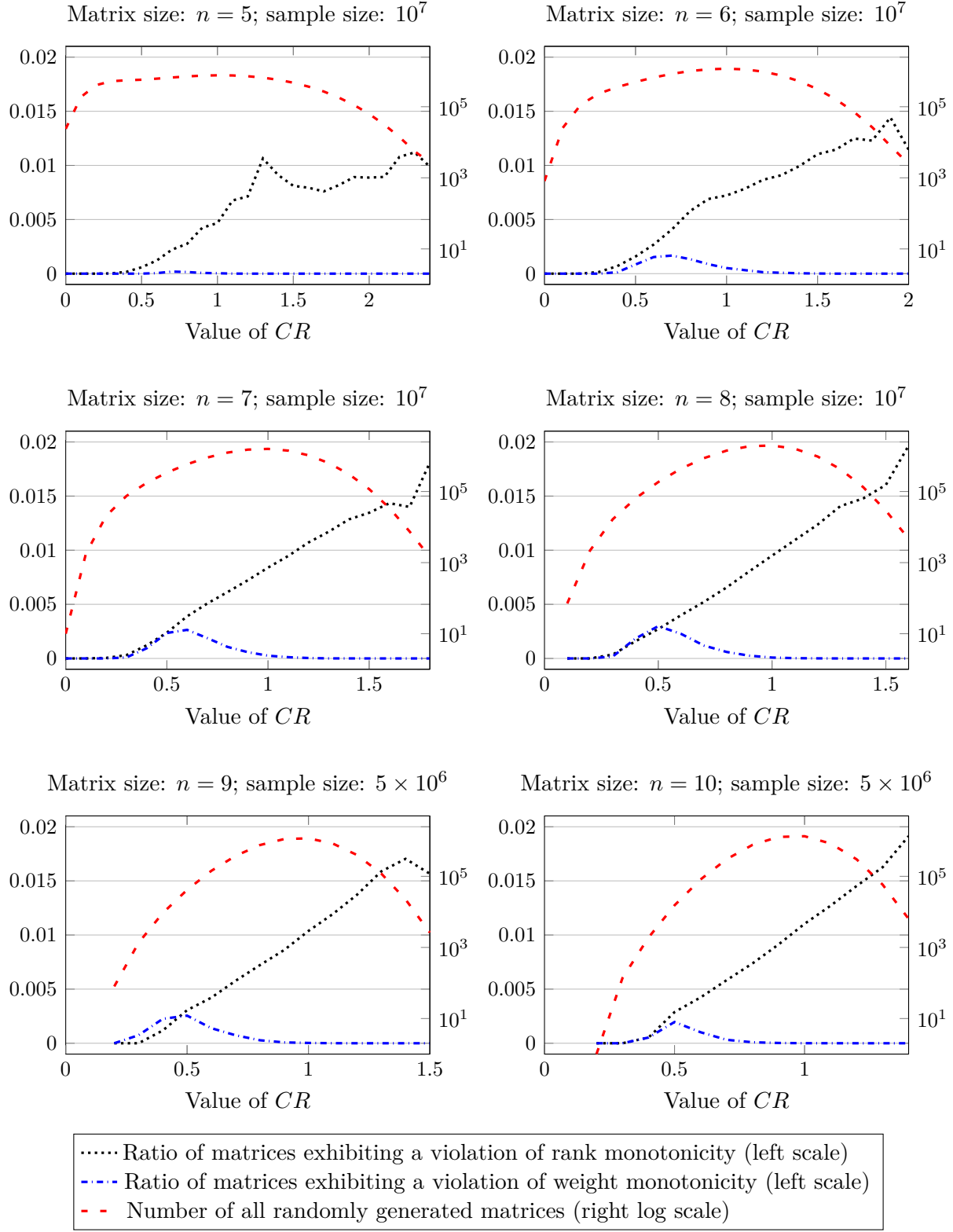
